\documentclass[12pt,leqno]{article}
\usepackage{latexsym}
\begin{titlepage}
\title{On the $p^{\lambda}$ problem }
\author{Stephan Baier
\thanks{e-mail: sbaier@mast.queensu.ca}}
\end{titlepage}
\begin{document}
\maketitle

\begin{center}
{\large \bf Abstract:}
\end{center}
We  deal with the distribution of the fractional parts of 
$p^{\lambda}$, $p$ running over the prime numbers and $\lambda$
being a fixed real number lying in the interval $(0,1)$.
Roughly speaking, we study the following question: Given a real $\theta$, how 
small may $\delta>0$ be choosen if we suppose that 
the number of primes $p\le N$ 
satisfying $\left\{p^{\lambda}-\theta\right\}<\delta$ is close to the 
expected one? 
We improve some results of Balog and Harman on this question 
for $\lambda<5/66$ if $\theta$ is rational and for $\lambda<1/5$ if 
$\theta$ is irrational. 

Our improvement is based on incorporating the zero 
detection argument into Harman$^{\prime}$s method and on using new mean 
value estimates for products of
shifted and ordinary (unshifted) Dirichlet polynomials.

\setcounter{page}{1}
\section{Introduction}
The well-known ``H conjecture'' (see [HaR]) states that $n^2+1$ is prime
infinitely many often. This is equivalent to the existence of
infinitely many primes $p$ satisfying $\{p^{1/2}\}<p^{-1/2}$. 
The current methods of analytic number theory are far from being sufficient to
prove these conjectures. 

However, Kubilius [Kub] and Ankeny
[Ank] proved already about fifty years ago that, 
assu\-ming the truth of the
Riemann Hypothesis for Hecke $L$-functions with Gr\"o\ss encharacters over 
${\bf Q}(i)$, 
$p=n^2+m^2$ is infinitely many often prime with $m\ll \log p$. This
implies that $\left\{p^{1/2}\right\}<p^{-1/2+\varepsilon}$ for infinitely many 
primes $p$. Of course, this is only a conditional result.
 
As demonstrated in [Ba1,2], [Ha1,2] and [BaH], 
it is also possible to obtain some {\bf unconditional} 
non-trivial results on small fractional parts of $p^{1/2}$, or more
generally, on small fractional parts of $p^{\lambda}$, $\lambda$ 
being a fixed real number lying in the interval $(0,1)$.
In particular, for $\lambda=1/2$ Balog and Harman
obtained $\left\{p^{1/2}\right\}<p^{-1/4+\varepsilon}$ for 
infinitely many primes $p$.
This result has recently be beaten. Combining Kubilius$^{\prime}$ 
ideas with efficient sieve 
methods, Harman and Lewis [HaL] unconditionally showed that the exponent 1/4 
may be replaced by 0.262. However, their method works only for $\lambda=1/2$,
whereas the methods in [Ba1,2], [Ha1,2] and [BaH] are applicable to all 
$\lambda$ in 
certain subintervals of $(0,1)$.   

In the present paper, we focus our interest mainly to small exponents
$\lambda$. Our starting point is the following 
result of Harman (Theorem 4 in
[Ha2]).\\
 
  {\sc Theorem 1:} \begin{it} Suppose that $\varepsilon>0$, 
$B>0$ and $\lambda\in (0,1/5]$ are given. Let $N\ge 3$. 
For every positive integer k define

  \begin{displaymath}
  e_1(\lambda,k):=\frac{5k-(2k+4)\lambda}{12k+4},
  \end{displaymath}

  \begin{displaymath}
  e_2(\lambda,k):=\frac{5k}{12k-6}-\lambda
  \end{displaymath}
and

  \begin{displaymath}
  e(\lambda,k):=\min\left\{e_1(\lambda,k),e_2(\lambda,k)\right\}.
  \end{displaymath}
Furthermore, define

  \begin{displaymath}
  E(\lambda):=\max\limits_{k\in {\bf N}}\ e(\lambda,k).
  \end{displaymath}
Then for

  \begin{equation}
  N^{-E(\lambda)+\varepsilon\lambda}\le \delta \le 1 \label{1}
  \end{equation}
we have

  \begin{equation}
  \sum\limits_{\scriptsize \begin{array}{cccc} &N<n\le 2N,&\\ 
  &\left\{n^{\lambda}\right\}<\delta&\end{array}} \hspace{-0.4cm}
  \Lambda(n)\ =\ \delta N \cdot\left(1+
  O\left(\frac{1}{(\log N)^B}\right)\right) \label{2}
  \end{equation}
as $N\rightarrow\infty$.\end{it}\\

Here, as in the following, $\Lambda(n)$ denotes the von Mangoldt
function.

As to be seen from the remark attached to Theorem 4 in [Ha2], this result
essentially keeps its validity if one introduces an additional summation
condition ``$\left[n^{\lambda}\right]\in {\bf A}$''
on the left side of (\ref{2}), where 
${\bf A}$ is any given subset of the 
set of positive integers (only the main term on the right side of (\ref{2}) 
correspondingly changes). Harman$^{\prime}$s motivation to introduce this 
additional condition appears to be the special case when ${\bf A}$ is the set 
of 
primes. 

Furthermore, in the same remark attached to Theorem 4 in [Ha2] it is
noted that the condition $\lambda\le 1/5$ may be replaced by $\lambda\le 1/2$ 
without any change in the result.  

To prove 
his result, Harman used density
estimates for the set of non-trivial zeta zeros and an estimate
for the $2k$-th power moment of Dirichlet polynomials 

\begin{displaymath}
\sum\limits_{m\sim M} a_m m^{it}.
\end{displaymath}

Hitherto, we have only considered small fractional parts of 
$p^{\lambda}$. 
A na\-tural genera\-lisation of this question is to consider small fractional parts
of $\left\{p^{\lambda}-\theta\right\}$, where $\theta$ is a given real number.
Unlike Theorem 4 in [Ha2], many results in [Ba1,2] and [Ha1,2] are formulated
for $\left\{p^{\lambda}-\theta\right\}$ with a general real $\theta$.
To extend Theorem 4 in [Ha2] in order to cover this general case, 
one needs estimates for 
power moments of {\bf shifted} 
Dirichlet polynomials 

$$
\sum\limits_{m\sim M} a_m(m+\theta)^{it}.
$$ 
In case
$\theta$ is rational these shifted Dirichlet polynomials 
can be easily rewritten as ordinary ones: If $\theta=b/q$, 
where 
$b$, $q$ are non-negative integers (without loss of generality, $\theta$ can
supposed to be non-negative), then

\begin{displaymath}
\sum\limits_{m\sim M} a_m(m+\theta)^{it} =q^{-it} \sum\limits_{m\sim M}
a_m(qm+b)^{it}.
\end{displaymath}
Therefore,
Harman$^{\prime}$s method works for all rational $\theta$, not only for
$\theta=0$. 

However, for irrational $\theta$ there
seem to be no reasonable known estimates of the
$2k$-th moment of $\sum\limits_{m\sim M} a_m(m+\theta)^{it}$ if $k>2$.  
Harman obtained such 
estimates only for $k\le 2$. But in this case his power moment 
estimates for irrational $\theta$ are essentially the same as the known ones 
for rational $\theta$. Thus, {\sc Theorem 1} keeps its validity also for
irrational $\theta$ if we replace the function $E(\lambda)$ by

  \begin{equation}
  E^*(\lambda):=\max \left\{e(\lambda,1),e(\lambda,2)\right\}. \label{3}
  \end{equation}

Summarising the above observations, 
{\sc Theorem 1} can be extended to the following\\

{\sc Theorem 2:} \begin{it}
Suppose that $\varepsilon>0$, $B>0$, $\lambda\in (0,1/2]$ 
and a real $\theta$ are given. Let $N\ge 3$.
Let ${\bf A}$ be an arbitrarily given subset of
the set of positive integers. Define $E(\lambda)$ as in {\sc Theorem} {\rm 1} 
and 
$E^*(\lambda)$ as in {\rm (\ref{3})}. 
Suppose that the condition {\rm (\ref{1})} is satisfied if $\theta$ is 
rational and that

  \begin{equation}
  N^{-E^*(\lambda)+\varepsilon\lambda}\le \delta \le 1 \label{4}
  \end{equation}
is satisfied if $\theta$ is irrational. Then we have

  \begin{equation}
  \sum\limits_{\scriptsize \begin{array}{cccc} &N<n\le 2N,&\\ 
  &\left\{n^{\lambda}-\theta\right\}<\delta,&\\ 
  &\left[n^{\lambda}\right]\in {\bf A}& \end{array}} \hspace{-0.4cm}
  \Lambda(n)\ =\ \frac{\delta}{\lambda} \cdot\hspace{-0.4cm}
  \sum\limits_{\scriptsize \begin{array}{cccc} 
  &N^{\lambda}<n\le (2N)^{\lambda},&\\ 
  &n\in {\bf A}& \end{array}} \hspace{-0.4cm}
  n^{1/\lambda-1}\ + \
  O\left(\frac{\delta N}{(\log N)^B}\right) \label{5}
  \end{equation}
as $N\rightarrow\infty$.\end{it}\\ 

It is easily verified that 
$E^*(\lambda)=5/14-2\lambda/7$ for $\lambda\le 5/18$. 
Using zero density estimates and trivially
estimating the shifted Dirichlet polynomials appearing in the method,
or directly applying Huxley$^{\prime}$s prime number theo\-rem,  
one obtains $5/12-\lambda$ in place of $E^*(\lambda)$, which yields a
better result
than $E^*(\lambda)$ if $\lambda<1/12$. This demonstrates that Harman$^{\prime}$s
method is ineffective if $\theta$ is irrational and $\lambda$ is
close to 0.

The first aim of the present paper is to prove a substantially better 
result than the one obtained from Huxley$^{\prime}$s prime number theorem for
irrational $\theta$ and $\lambda$ close to 0. 
Our second aim is to improve {\sc Theorem 2} for 
rational $\theta$. We shall prove the following \\

{\sc Theorem 3:} \begin{it} 
Suppose that $\varepsilon>0$, $B>0$, $\lambda\in (0,1/2]$ 
and a real $\theta$ are given. If $\theta$ is irrational, then suppose 
that $\lambda< 5/19$. Let $N\ge 3$. 
Let ${\bf A}$ be an arbitrarily given subset of
the set of positive integers. For every positive integer k define

  \begin{displaymath}
  f_1(\lambda,k):=\frac{5}{12}-\frac{k+6}{6(k+1)}\cdot\lambda,
  \end{displaymath}

  \begin{displaymath}
  f_2(\lambda,k):=\frac{5}{11}-\frac{5k+1}{11}\cdot \lambda
  \end{displaymath}
and 

  \begin{displaymath}
  f(\lambda,k):=\min\left\{f_1(\lambda,k),f_2(\lambda,k)\right\}.
  \end{displaymath}
Furthermore, define

  \begin{displaymath}
  F_{\theta}(\lambda):=\left\{ \begin{array}{llll} 
  F(\lambda) & \mbox{ if } \theta \mbox{ is 
  rational,} \\ \\ f(\lambda,1) & \mbox{ otherwise,} \end{array} \right.
  \end{displaymath}
where 

  \begin{displaymath}
  F(\lambda):=\max\limits_{k\in {\bf N}}\ f(\lambda,k).
  \end{displaymath}
Suppose that

  \begin{equation}
  N^{-F_{\theta}(\lambda)+\varepsilon\lambda}\le \delta \le 1. \label{6}
  \end{equation}
Then we have the asymptotic estimate {\rm (5)} as 
$N\rightarrow \infty$.\end{it}\\

We note that 

\begin{equation}
f(\lambda,1)=f_1(\lambda,1)=\frac{5}{12}-\frac{7\lambda}{12} 
\end{equation}
for all $\lambda>0$.\\ 

In the next section we shall discuss {\sc Theorem 3} in detail and compare
this result with {\sc Theorem 2}. From the third section onwards we shall 
prove {\sc Theorem 3}.\\

{\bf Acknowledgements.} The main part of research was carried out when the 
author held a
Postdoctoral Fellowship at the Harish-Chandra Research Institute at
Allahabad (India). 
The paper was completed when he held a Marie-Curie Postdoctoral
Fellowship at the Department of Pure Mathematics and Mathematical Statistics 
of the University of Cambridge (United Kingdom). 
He thanks the Harish-Chandra Research Institute as well as the Marie-Curie Foundation of the European Union for financial 
support and the two above-mentioned institutes for providing very good
working conditions. Furthermore, he wishes to thank the Referee for some 
useful comments as well as Glyn Harman for bringing the paper [HaL] to the
author$^{\prime}$s attention. 

\section{Discussion of Theorem 3}

In case $\theta$ is irrational it is supposed in {\sc Theorem 3} that 
$\lambda< 5/19$. We have 
$E^*(\lambda)=5/14-2\lambda/7$ if $\lambda< 5/19$. Thus, we get 
$5/12-7\lambda/12=f(\lambda,1)>E^*(\lambda)$ 
if $\lambda<1/5$. Therefore, 
{\sc Theorem 3}  yields a sharper result than {\sc Theorem 2} if $\theta$ is 
irrational and $\lambda<1/5$. 
Moreover, since $f(\lambda,1)>5/12-\lambda$ for all $\lambda>0$, 
{\sc Theorem 3} is
always sharper than the result directly 
obtained from Huxley$^{\prime}$s prime number 
theorem (see the preceding  section).

We now turn to the case when $\theta$ is rational. It is easily 
verified that
$F(\lambda)=5/12-\lambda/6+O\left(\lambda^2\right)$ as
$\lambda\rightarrow 0$, whereas
$E(\lambda)=5/12-\lambda/2+O\left(\lambda^2\right)$ as
$\lambda\rightarrow 0$.  
Thus, $F(\lambda)>E(\lambda)$ for
every sufficiently small $\lambda>0$. Therefore, {\sc Theorem 3} is sharper 
than {\sc Theorem 2} if $\theta$ is rational and $\lambda$ is sufficiently 
small. We now make this observation more precise by analysing and comparing 
$E(\lambda)$ and $F(\lambda)$.

It is easily seen that there are sequences $\left(\eta_{1,k}\right)$, 
$\left(\eta_{2,k}\right)$ of real numbers with

$$
1/2=\eta_{2,1}>\eta_{1,2}>
\eta_{2,2}>\eta_{1,3}>\eta_{2,3}>...
$$ 
and 

$$
\lim\limits_{k\rightarrow\infty} \eta_{i,k}=0\ \ \ \ (i=1,2),
$$
such that

$$
E(\lambda)=\left\{\begin{array}{llll} e_1(\lambda,k-1) 
& \mbox{ if } \eta_{2,k-1}\ge
\lambda\ge \eta_{1,k} \\ \\  e_2(\lambda,k) & \mbox{ if } \eta_{1,k}\ge
\lambda\ge \eta_{2,k}, \end{array} \right.
$$
where the functions $e_i(\lambda,k)$ $(i=1,2)$ are defined as in 
{\sc Theorem 1}. So 
$E(\lambda)$ is a continuous piecewise linear function. 
To determine
$\eta_{i,k}$ for $k\ge 2$, we simply have to solve the linear equations

$$
e_1(x,k-1)=e_2(x,k) 
$$  
and 

$$
e_2(x,k)=e_1(x,k). 
$$  
In this manner, we obtain 

$$
\eta_{1,k}=\left(\frac{5}{12}+\frac{1}{6(k-1)}\right)\cdot \frac{1}{k-1/2}
$$
and

$$
\eta_{2,k}=\frac{5}{12}\cdot \frac{1}{k-1/2}
$$
if $k\ge 2$. Similarly, we obtain 

\begin{equation}
F(\lambda)=\left\{\begin{array}{llll} f_1(\lambda,k-1) & \mbox{ if } 
\phi_{2,k-1}\ge
\lambda\ge \phi_{1,k}, \\ \\  f_2(\lambda,k) & \mbox{ if } \phi_{1,k}\ge
\lambda\ge \phi_{2,k}, \end{array} \right. \label{8}
\end{equation}
where $\phi_{2,1}=1/2$ and 

$$
\phi_{1,k}=\frac{1}{2(6k-1-11/k)}
$$
and

$$
\phi_{2,k}=\frac{1}{2(6k-1-11/(k+1))}
$$
if $k\ge 2$.
Using these explicit expressions for $E(\lambda)$ and
$F(\lambda)$, it is not difficult 
to calculate that $F(\lambda)>E(\lambda)$ whenever $\lambda<5/66$ or 
$1/3<\lambda<1/2$.
Consequently, for rational $\theta$
{\sc Theorem 3}  yields a sharper result than {\sc Theo\-rem 2} in these 
$\lambda$-ranges. \newpage

{\bf Some additional remarks.} \\ \\
a) We have $E^*(1/2)=E(1/2)=f(1/2,1)=F(1/2)=1/8$ and 
$E(\lambda)=F(\lambda)=5/11-\lambda$ if $5/66=\eta_{2,6}\le \lambda 
\le 9/110=\eta_{1,6}$.\medskip\\
b) In the homogeneous 
case ${\bf A}={\bf N}$, Harman and Balog were able to prove 
better results than {\sc Theorems 2} and 3 for  $\lambda>1/5$ (see Theorem
3 of [Ha2], Theorem 2 of [BaH] and the papers [Ba1,2], [Ha1]). \medskip\\ 
c) It is not difficult to prove that on the Riemann Hypothesis
the function $F_{\theta}(\lambda)$ in {\sc Theorem 3} can be replaced by
$(1-\lambda)/2$ for all $\lambda$ in the interval $0<\lambda<1$ and all real 
$\theta$.

\section{Auxiliary results and outline of the method}
  
We let
the conditions of {\sc Theorem 3} be kept throughout the remaining part of 
the paper. Without loss of generality,
we continually suppose that $0\le \theta<1$. By $\varepsilon$ and $B$
we always mean 
the constants $\varepsilon$ and $B$ from {\sc Theorem 3}.  
 
We define $h=h_{\theta}(\lambda)$ to be the smallest
positive integer such that

\begin{equation}
F_{\theta}(\lambda)=f(\lambda,h). \label{102}
\end{equation}
In particular, we always have $h=1$ when $\theta$ is irrational. When $\theta$
is rational, we have $h\ge 2$ if and only if $\lambda<1/11=\phi_{1,2}$. 
In this case, we obtain 

\begin{equation}
h\lambda< 2/11 \label{9}
\end{equation} 
from 
(\ref{8}) by a short calculation. 

Using  (\ref{102}) and
the definition of $F_{\theta}(\lambda)$ in {\sc Theorem 3}, 
we obtain 

  \begin{equation}
  \frac{5}{12}\le \lambda+f(\lambda,h)\le \frac{5}{11} \label{10} 
  \end{equation} 
if $h\ge 2$ and  

  \begin{equation}
  \frac{\lambda}{2}+f(\lambda,h)\le \frac{5}{11} \label{104}
  \end{equation}
in any case. Moreover, we have 

  \begin{equation}
  \lambda+f(\lambda,1)=\lambda+f_1(\lambda,1)=
  \frac{5(1+\lambda)}{12}\le \frac{5}{8} \label{103}
  \end{equation}
for every $\lambda\le 1/2$ and 

  \begin{equation}
  \lambda\le \frac{\lambda+f(\lambda,1)-\varepsilon\lambda}{2} \label{105}
  \end{equation} 
if $\lambda<5/19$ and $\varepsilon \le (5-19\lambda)/(12\lambda)$. We shall use the 
inequalities (10)-(14) in the course of this paper. 

Next, we introduce some more notations.
We write
  \begin{displaymath}
  D_y(u,s):=\hspace{-0.4cm}\sum\limits_{\scriptsize \begin{array}{cccc}
  & N^{\lambda}<n\le u,& \\
  &n\in {\bf A}& \end{array}}\hspace{-0.4cm}
  (n+y)^{s/\lambda-1}
  \end{displaymath}
for any real $u> N^{\lambda}$, $y\ge 0$ and complex $s$. Moreover, we put
  
  \begin{displaymath}
  D_y(s):=D_y((2N)^{\lambda},s).
  \end{displaymath}
As usual, by the
symbol $\rho$ we denote the non-trivial zeta zeros, and we write
$\gamma$ for the imaginary part and $\beta$ for the real part of
$\rho$. We define

  \begin{equation}
  S_{\theta}(u,\sigma):=
  \hspace{-0.4cm} \sum\limits_{\scriptsize \begin{array}{cccc} 
  &\rho:\ \! 0<\gamma\le T,&\\ &\sigma \le \beta\le \sigma+1/(\log N)& 
  \end{array}}\hspace{-0.4cm}  
  \left\vert D_{\theta}(u,i\gamma) \right\vert \label{11}
  \end{equation}
for any $\sigma$ with $0\le \sigma\le 1-1/(\log N)$, where the parameter $T$ 
shall be fixed
at the beginning of the next section.  

The first step of our method is to use 
the explicit formula of Landau in order to reduce the sum on the left side of
(\ref{5}) to sums of the form $S_{\theta}(u,\sigma)$.\\

  {\sc Proposition 1:} (explicit formula) \begin{it}
  For $x>2$, $T_0>1$ we have

  \begin{displaymath}
  \sum\limits_{n < x} \Lambda(n)= x-
  \sum\limits_{\rho:\ \! \vert \gamma \vert \le T_0} \frac{x^{\rho}}{\rho}
  + O\left(\frac{x}{T_0}\cdot  (\log xT_0)^2+ \log x\right).
  \end{displaymath} \end{it}\medskip

We then estimate the sum $S_{\theta}(u,\sigma)$ in several $\sigma$-regions
by different me\-thods.   

To control the range $0\le \sigma\le 6/11$, we  use the following mean
value estimate for shifted
Dirichlet polynomials which can be established
in the same manner as the corresponding well-known mean value estimate
for ordinary (unshifted) Dirichlet polynomials (see [Ivi] for example).\\

{\sc Proposition 2:} \begin{it} 
Suppose that  $0\le \theta <1$, $K\ge 1$ and  $T\ge 1$. 
Let  $(a_{k})$ be an arbitrary sequence of complex
numbers. Suppose that $\vert a_k \vert \le A$ for all $k\sim K$. 
Let $(t_r)$ be a
monotone increasing sequence of positive real numbers, such that 
$t_{r+1}-t_r\ge 1$ for every positive integer $r$. Let $R$ be
a positive integer. Suppose that $t_R\le T$. Then,

  \begin{displaymath}
  \sum\limits_{r=1}^R \left\vert \sum\limits_{k\sim K} 
  a_{k}(k+\theta)^{it_r} \right\vert^2 \ \ll \
  A^2(T+K)K\log(2K). 
  \end{displaymath}\end{it}

Here, as in the following, the notation $k \sim K$ means $K< k\le 2K$.

To tackle $S_{\theta}(u,\sigma)$ in the range $1-\Delta<\sigma\le 1$, $\Delta$
being defined as in (\ref{500}), 
we use 
the
second zero density estimate of the later {\sc Proposition 6} as well as
Vinogradov$^{\prime}$s zero-free-region result.\\

{\sc Proposition 3 (Vinogradov):} (see [Ivi]) \begin{it}
There is an absolute constant $C>0$ such that

  \begin{displaymath}
  \beta\le 1-
  C(\log \vert \gamma \vert)^{-2/3}(\log \log \vert \gamma \vert)^{-1/3}
  \end{displaymath}
for every non-trivial zeta zero $\rho=\beta+i\gamma$.
\end{it}\\

To calculate $S_{\theta}(u,\sigma)$ in the range 
$6/11< \sigma \le 1-\Delta$, we employ the 
following relation which is also the basis of
the zero detection method for counting non-trivial zeta zeros 
(cf. [Ivi] for example).\\

{\sc Proposition 4:} \begin{it} Suppose that $X,Y\ge 1$,
  $T>1$, 
$\log N \ll \log T \ll \log Y \ll \log T \ll \log N$ and $\log X \ll
  \log N$. Define 
  
  \begin{displaymath}
  M_X(s):=\sum\limits_{n\le X} \mu(n)n^{-s}
  \end{displaymath}
and

  \begin{displaymath}
  a(k):=\hspace{-0.4cm}
  \sum\limits_{\scriptsize \begin{array}{cccc} &d\vert k& \\
  &d\le X& \end{array}} \hspace{-0.4cm} \mu(d).
  \end{displaymath}
Then we have $a(k)=0$ if $1<k\le X$,  

  \begin{displaymath}
  \zeta(s)M_X(s)=\sum\limits_{k=1}^{\infty} a(k)k^{-s}
  \end{displaymath}
if {\rm Re} $s>1$ and 

  \begin{equation}
  \frac{1}{2} < U_1(\rho) \hspace{0.5cm} \mbox{or} \hspace{0.5cm} \frac{1}{2}
  < U_2(\rho) \label{12}
  \end{equation}\\
if $N$ is sufficiently large and $\rho=\beta+i\gamma$ is a non-trivial
zeta zero satisfying 
$\beta\ge 1/2$ and $(\log N)^2<\gamma\le T$, where
  \begin{eqnarray}
  U_1(\rho) &:=& \left\vert \ \int\limits_{-(\log N)^2}^{(\log N)^2}
  \zeta(1/2+i(\gamma+t))M_X(1/2+i(\gamma+t))\cdot \right. \label{13} \\ & &  
  \left. \begin{array}{llll} \\ Y^{1/2-\beta+it} 
  \Gamma(1/2-\beta+it)\ {\rm d}t \\ \\ \end{array} \right\vert  \nonumber
  \end{eqnarray}
and

  \begin{equation}
  U_2(\rho) :=\left\vert \sum\limits_{X<k\le Y(\log N)^2} 
  a(k)k^{-\rho}e^{-k/Y} \right\vert. \label{14}
  \end{equation}
The sum on the right side of {\rm (\ref{14})}
is supposed to equal 0 if $X\ge Y(\log N)^2$.
We call every zero $\rho$ satisfying $1/2< U_j(\rho)$ a zero of
type $j$ $(j=1,2)$. 
\end{it}\\

We shall now similarly proceed as in the original zero-detection
method, with the diffe\-rence that here every non-trivial zeta
zero $\rho$ is 
weighted
by $\vert D_{\theta}(u,i\gamma)\vert$. Our method shall lead
us to the problem of estimating mean values 
of products of shifted and ordinary (unshifted) Dirichlet polynomials. In the
following, we state such mean value estimates.\\

{\sc Theorem 4:} \begin{it} 
Suppose that  $\alpha\not=0$, $0\le \theta<1$, $T> 0$, 
$K\ge 1$, $L\ge 1$. If $\theta\not=0$, then additionally suppose that 
$L\le T^{1/2}$. 
Let  $(a_{k})$ and  $(b_{l})$ be arbitrary sequences of complex
numbers. Suppose that $\vert a_k \vert \le A$ for all $k\sim K$ and 
$\vert b_l \vert \le B$ for all $l\sim L$. Then,
  \begin{eqnarray}
  & & \int\limits_0^T \left\vert \sum\limits_{k\sim K} 
  a_{k}k^{it} \right\vert^2 \
  \left\vert \sum\limits_{l\sim L} 
  b_{l}(l+\theta)^{i\alpha t} \right\vert^2  {\rm d}t \label{15} \\ 
  \nonumber\\
  &\ll& A^2B^2(T+KL)KL\log^3(2KLT), \nonumber
  \end{eqnarray}
the implied $\ll$-constant depending only on $\alpha$. If $\theta=0$, then
$\log^3(2KLT)$ on the right side of {\rm(\ref{15})} may be replaced by 
$\log^2(2KLT)$.\end{it}\\

We shall rather need a discrete form of {\sc Theorem 4}, namely\\

{\sc Theorem $4^*$:} \begin{it} 
Let the conditions of {\sc Theorem 4} be kept. Moreover, let $(t_r)$ be a
monotone increasing sequence of positive real numbers, such that 
$t_{r+1}-t_r\ge 1$ for every positive integer $r$. Further, let $R$ be
a positive integer. Suppose that $t_R\le T$. Then,
  \begin{eqnarray}
  & & \sum\limits_{r=1}^R \left\vert \sum\limits_{k\sim K} 
  a_{k}k^{it_r} \right\vert^2 \
  \left\vert \sum\limits_{l\sim L} 
  b_{l}(l+\theta)^{i\alpha t_r} \right\vert^2 \label{16} \\ \nonumber\\
  &\ll& A^2B^2(T+KL)KL\log^4(2KLT), \nonumber
  \end{eqnarray}
the implied $\ll$-constant depending only on $\alpha$. If $\theta=0$, then
$\log^4(2KLT)$ on the right side of {\rm (\ref{16})} may be replaced by 
$\log^3(2KLT)$.\end{it}\\

We postpone the proofs of {\sc Theorems 4}, $4^*$ to the last section, in 
which we shall also derive the following more general mean value estimate 
from {\sc Theorem} $4^*$.\\

{\sc Theorem 5:} \begin{it} 
Suppose that  $\alpha\not=0$, $0\le \theta<1$, $T> 0$, 
$1\le K_1<K_2$, $1\le L_1<L_2$, $h\in {\bf N}$, $\varepsilon_0>0$ and
$1/2\le \sigma\le 1$. If $\theta$ is irrational, then additionally suppose 
that $L_2\le T^{1/2}$ and $h=1$. 
Let  $(a_{k})$ and  $(b_{l})$ be arbitrary sequences of complex
numbers. Suppose that $\vert a_k \vert \le A$ and 
$\vert b_l \vert \le B$ for all positive integers $k$ and $l$. Let $(t_r)$ 
be a monotone increasing sequence of positive real numbers, such that 
$t_{r+1}-t_r\ge 1$ for every positive integer $r$. Let $R$ be
a positive integer. Suppose that $t_R\le T$. Then we have
  \begin{eqnarray}
  & & \sum\limits_{r=1}^R \left\vert \sum\limits_{K_1<k\le K_2} 
  a_{k}k^{-(\sigma+it_r)} \right\vert^2 \
  \left\vert \sum\limits_{L_1<l\le L_2} 
  b_{l}(l+\theta)^{i\alpha t_r-1} \right\vert^{2h}  \label{17} 
  \\ \nonumber\\
  &\ll&
  A^2B^{2h}\left(TK_1^{1-2\sigma}L_1^{-h}+K_2^{2(1-\sigma)}\right)
  (K_2L_2T)^{\varepsilon_0}, \nonumber
  \end{eqnarray}
the implied $\ll$-constant depending only on $\alpha$, $\theta$, $h$ and 
$\varepsilon_0$.\end{it}\\

{\sc Theorem 5} is made for a direct application in the present paper. 

In
addition to these mean value estimates, we use the following
well-known fourth power moment estimate for the Riemann zeta function
on the critical line.\\

{\sc Proposition 5:} We have

  \begin{displaymath}
  \int\limits_0^{T} \vert \zeta(1/2+it) \vert^4\ {\rm d}t\ \ll \ T\log^4 T.
  \end{displaymath}

Finally, we shall employ zero density estimates of Ingham and Huxley 
which themselves are consequences of the zero-detection
method.\\ 

{\sc Proposition 6:} (see [Ivi]) \begin{it} For $T>2$ we have

  \begin{displaymath}
  {\bf N}(\sigma,T)\ll \left\{ \begin{array}{llll} T^{3(1-\sigma)/(2-\sigma)}(\log T)^{5} &
  \mbox{\rm if }\ 1/2\le \sigma \le 3/4, \\ \\
  T^{3(1-\sigma)/(3\sigma-1)}(\log T)^{44} & \mbox{\rm if }\ 3/4\le
  \sigma \le 1, \end{array} \right.
  \end{displaymath}
where ${\bf N}(\sigma,T)$ denotes the number of zeta zeros
$\rho=\beta+i\gamma$ with $\beta\ge \sigma$ and $0<\gamma\le T$.
\end{it}
      
\section{Reduction to sums over nontrivial zeta zeros}

We define 

\begin{displaymath}
  T_0:= \frac{N^{\lambda}(\log N)^{B+2}}{\delta} 
  \end{displaymath}
and 
 
  \begin{equation}
  T:= N^{\lambda+f(\lambda,h)-\varepsilon\lambda}(\log N)^{B+2}.
  \label{18}
  \end{equation}
We note that $T_0\le T$ by (\ref{102}) and 
condition (\ref{6}) of {\sc Theorem 3}. Furthermore,
we state the following five bounds, which shall be used in the course of 
this paper. If $h\ge 2$, then we have 

  \begin{equation}
  T\ll N^{5/11-\varepsilon\lambda/2} \label{100}
  \end{equation}
as well as 

  \begin{equation}
  N^{5/12-\varepsilon\lambda} \ll T \label{101}
  \end{equation}
by (\ref{10}). We always have

  \begin{equation}
  TN^{-\lambda/2}\ll N^{5/11-\varepsilon\lambda/2} \label{107}
  \end{equation}
by (\ref{104}). If $h=1$, then we have 
  \begin{eqnarray}
  T&=&N^{5(1+\lambda)/12-\varepsilon\lambda}(\log N)^{B+2} \label{106}\\
  &\ll& N^{5/8-\varepsilon\lambda/2}\nonumber
  \end{eqnarray}
by (\ref{103}).  If $h=1$, $\lambda<5/19$, 
$\varepsilon \le (5-19\lambda)/(12\lambda)$ and $N\ge 5$, then we have

  \begin{equation}
  (2N)^{\lambda}\le T^{1/2}\label{109}
  \end{equation} 
by (\ref{105}). 

By means of {\sc Proposition 1}, we now decompose the sum on the
left side of (\ref{5}) 
into a main term and an error term involving non-trivial
zeta zeros.\\
 
{\sc Lemma 1:} \begin{it} We have
  \begin{eqnarray*}
  & & \hspace{-0.4cm}
  \sum\limits_{\scriptsize \begin{array}{cccc} &N<n\le 2N,&\\ 
  &\left\{n^{\lambda}-\theta\right\}<\delta,&\\ 
  &\left[n^{\lambda}\right]\in {\bf A}& \end{array}} \hspace{-0.4cm}
  \Lambda(n)\  -\  \frac{\delta}{\lambda} \cdot D_0(1) 
  \\ &\ll& \delta (\log N)^{B+3} \sup\limits_{0\le \sigma \le 1-1/(\log N)}\
  N^{\sigma}\ \sup\limits_{N^{\lambda}<u \le (2N)^{\lambda}} 
  S_{\theta}(u,\sigma) \\ & &  
  +\ \frac{\delta N}{(\log N)^B}\ +\  N^{\lambda}\log N,
  \end{eqnarray*}
the implied $\ll$-constant depending only on $\lambda$ and $B$.
\end{it}\\ 

{\bf Proof:} 
Obviously, the sum in question can be written in the form
  \begin{eqnarray}
  \sum\limits_{\scriptsize \begin{array}{cccc} &N<n\le 2N,&\\ 
  &\left\{n^{\lambda}-\theta\right\}<\delta,&\\ 
  &\left[n^{\lambda}\right]\in {\bf A}& \end{array}} \hspace{-0.4cm}
  \Lambda(n) &=& \hspace{-0.4cm}
  \sum\limits_{\scriptsize \begin{array}{cccc} 
  &N^{\lambda}<n\le (2N)^{\lambda},& \\ 
  &n\in {\bf A}& \end{array}}\hspace{-0.4cm}
  \ \sum\limits_{(n+\theta)^{1/\lambda}\le m< 
  (n+\theta+\delta)^{1/\lambda}} \Lambda(m)\
  \label{19} \\ & &  \hspace{0.1cm} 
  +\ O\left(N^{\lambda}\log N\right).\nonumber
  \end{eqnarray}
Combining this estimate and {\sc Proposition 1}, and taking the condition 
(\ref{6})
into account, we get
  \begin{eqnarray}
  & & \sum\limits_{\scriptsize \begin{array}{cccc} &N<n\le 2N,&\\ 
  &\left\{n^{\lambda}-\theta\right\}<\delta,&\\ 
  &\left[n^{\lambda}\right]\in {\bf A}& \end{array}} \hspace{-0.4cm}
  \Lambda(n) \label{20}\\ \nonumber\\
  &=&\hspace{-0.4cm}\sum\limits_{\scriptsize \begin{array}{cccc}
  & N^{\lambda}<n\le (2N)^{\lambda},& \\
  &n\in {\bf A}& \end{array}}\hspace{-0.4cm}
  \left((n+\theta+\delta)^{1/\lambda}-(n+\theta)^{1/\lambda}\right)
  \nonumber\\ & & -\hspace{-0.4cm}\sum\limits_{\scriptsize \begin{array}{cccc}
  & N^{\lambda}<n\le (2N)^{\lambda},& \\ 
  &n\in {\bf A}& \end{array}}\hspace{-0.4cm}
  \sum\limits_{\rho:\ \! \vert \gamma \vert \le T_0}
  \frac{(n+\theta+\delta)^{\rho/\lambda}-
  (n+\theta)^{\rho/\lambda}}{\rho} \nonumber\\ \nonumber\\& &
  +\ O\left(\frac{\delta N}{(\log N)^B}+N^{\lambda}\log N\right).\nonumber
  \end{eqnarray}
Using Taylor$^{\prime}$s formula,
we approximate the first sum on the right side of (\ref{20}) by
  \begin{eqnarray}\label{21} \\
  \hspace{-0.4cm}\sum\limits_{\scriptsize \begin{array}{cccc}
  & N^{\lambda}<n\le (2N)^{\lambda},& \\
  &n\in {\bf A}& \end{array}}\hspace{-0.4cm}
  \left((n+\theta+\delta)^{1/\lambda}-(n+\theta)^{1/\lambda}\right)
  =
  \frac{\delta}{\lambda} \cdot D_0(1)
  \ +\  O\left(\delta N^{1-\lambda}\right).\nonumber
  \end{eqnarray}
The fraction within the double
sum on the right side of (\ref{20}) can be written as an integral, namely

  \begin{displaymath}
  \frac{(n+\theta+\delta)^{\rho/\lambda}-
  (n+\theta)^{\rho/\lambda}}{\rho}=
  \frac{1}{\lambda} \cdot \int\limits_{\theta}^{\theta+\delta} 
  (n+y)^{\rho/\lambda-1}\ {\rm d}y.
  \end{displaymath}
From that and  the symmetry of the set of zeta zeros 
it follows that the double sum on the right side of (\ref{20}) 
can be estimated by
  \begin{eqnarray}
  & & \hspace{-0.4cm}\sum\limits_{\scriptsize \begin{array}{cccc}
  & N^{\lambda}<n\le (2N)^{\lambda},& \\
  &n\in {\bf A}& \end{array}}\hspace{-0.4cm}
  \sum\limits_{\rho:\ \! \vert \gamma \vert \le T_0}
  \frac{(n+\theta+\delta)^{\rho/\lambda}-
  (n+\theta)^{\rho/\lambda}}{\rho} \label{22} \\ \nonumber\\ &\ll&
  \delta \ \sup\limits_{0\le y\le \delta}\
  \sum\limits_{\rho:\ \! 0< \gamma \le T_0} 
  \left\vert D_{\theta+y}\left(\rho\right)\right\vert.\nonumber
  \end{eqnarray}
Moreover, we have
  \begin{eqnarray}
  & & \sum\limits_{\rho:\ \! 0< \gamma \le T_0} 
  \left\vert D_{\theta+y}\left(\rho\right)\right\vert \label{23}\\ \nonumber\\ 
  &\ll& (\log N)\sup\limits_{0\le \sigma\le 1-1/(\log N)} \hspace{-0.4cm}
  \sum\limits_{\scriptsize \begin{array}{cccc} 
  &\rho:\ \! 0<\gamma\le T_0,&\\ &\sigma \le \beta\le \sigma+1/(\log N)& 
  \end{array}}\hspace{-0.4cm}  
  \left\vert D_{\theta+y}\left(\rho\right) \right\vert.\nonumber
  \end{eqnarray}

The next step is to reduce the shifted Dirichlet polynomial
$D_{\theta+y}(\rho)$ to $D_{\theta}(i\gamma)$. By partial summation, we get
  \begin{eqnarray}
  & & \label{24} \\
  D_{\theta+y}(\rho) &=&
  \hspace{-0.4cm}\sum\limits_{\scriptsize \begin{array}{cccc}
  & N^{\lambda}<n\le (2N)^{\lambda},& \\
  &n\in {\bf A}& \end{array}}\hspace{-0.4cm} (n+\theta+y)^{\beta/\lambda}
  \left(1+\frac{y}{n+\theta}\right)^{i\gamma/\lambda-1}
  (n+\theta)^{i\gamma/\lambda-1}
  \nonumber\\ &=&
  \left((2N)^{\lambda}+\theta+y\right)^{\beta/\lambda}
  \left(1+\frac{y}{(2N)^{\lambda}+\theta}\right)^{i\gamma/\lambda-1}
  D_{\theta}(i\gamma)
  \nonumber\\ & &-
  \int\limits_{N^{\lambda}}^{(2N)^{\lambda}} \frac{{\rm d}}{{\rm d}u}
  \left((u+\theta+y)^{\beta/\lambda}
  \left(1+\frac{y}{u+\theta}\right)^{i\gamma/\lambda-1}\right) 
  D_{\theta}(u,i\gamma)
  \ {\rm d}u.\nonumber
  \end{eqnarray}
For $N^{\lambda}\le u \le (2N)^{\lambda}$, $0\le y\le \delta$,
$0<\gamma\le T_0$ and $0\le \sigma\le \beta\le \sigma+1/(\log N)\le 1$ we 
have
  \begin{eqnarray}
  & & \frac{{\rm d}}{{\rm d}u}
  \left((u+\theta+y)^{\beta/\lambda}
  \left(1+\frac{y}{u+\theta}\right)^{i\gamma/\lambda-1}\right) \label{25} 
  \\
  &=& \frac{\beta}{\lambda}\cdot (u+\theta+y)^{\beta/\lambda-1}
  \left(1+\frac{y}{u+\theta}\right)^{i\gamma/\lambda-1}-\nonumber\\ & &
  \left(\frac{i\gamma}{\lambda}-1\right) \cdot \frac{y}{(u+\theta)^2}
  \cdot (u+\theta+y)^{\beta/\lambda}
  \left(1+\frac{y}{u+\theta}\right)^{i\gamma/\lambda-2} \nonumber\\ \nonumber\\
  &\ll& N^{\beta-\lambda}+\delta T_0 N^{\beta-2\lambda}\nonumber\\ 
  &\ll& N^{\sigma-\lambda}(\log N)^{B+2}.\nonumber
  \end{eqnarray}
From (\ref{24}) and (\ref{25}), we obtain
  \begin{eqnarray}
  & & \label{26} \\ \vert D_{\theta+y}(\rho)\vert 
  &\ll& N^{\sigma}(\log N)^{B+2}\cdot
  \left(\left\vert D_{\theta}(i\gamma)\right\vert+ N^{-\lambda} \cdot
  \int\limits_{N^{\lambda}}^{(2N)^{\lambda}}
  \left\vert D_{\theta}(u,i\gamma) \right\vert
  {\rm d}u\right).\nonumber
  \end{eqnarray}
From (\ref{26}), 
$T_0\le T$ and the definition of 
$S_{\theta}(u,\sigma)$ in (\ref{11}), we derive
  \begin{equation}
  \hspace{-0.4cm} \sum\limits_{\scriptsize \begin{array}{cccc} 
  &\rho:\ \! 0<\gamma\le T_0,&\\ &\sigma \le \beta\le \sigma+1/(\log N)& 
  \end{array}}\hspace{-0.4cm}  
  \left\vert D_{\theta+y}\left(\rho\right) \right\vert \ \ll \
   N^{\sigma}(\log N)^{B+2}\cdot
  \sup\limits_{N^{\lambda}<u\le (2N)^{\lambda}}
  S_{\theta}(u,\sigma).\label{27}
  \end{equation}

Combining (\ref{19}), (\ref{20}), (\ref{21}), (\ref{22}), (\ref{23}) 
and (\ref{27}), we obtain the
desired estimate. $\Box$\\

By {\sc Lemma 1}, in order to prove (\ref{5}), we still have to show that

  \begin{equation}
  N^{\sigma}S_{\theta}(u,\sigma)\ll \frac{N}{(\log N)^{2B+3}} \label{28}
  \end{equation}
for $0\le \sigma\le 1-1/(\log N)$. This shall be the task of the next sections.
$\Box$

\section{Estimation of $N^{\sigma}S_{\theta}(u,\sigma)$ for  
$0\le\sigma\le 6/11$ and for $1-\Delta<\sigma\le 1$} 

In this section we establish (\ref{28}) for $0\le\sigma\le 6/11$ and for
$1-\Delta<\sigma\le 1$, where

\begin{equation}
\Delta:=\left\{\begin{array}{llll} 0.212 & \mbox{ if } h\ge 2, \\ \\
1/36 & \mbox{ if } h=1. \end{array} \label{500}\right.
\end{equation}

{\sc Lemma 2:} \begin{it} Without loss of generality assume that 
$\varepsilon\le f(\lambda,h)/\lambda$. Then for $0\le \sigma\le 6/11$ we have

  \begin{displaymath}
  N^{\sigma}S_{\theta}(u,\sigma) \ll N^{1-\varepsilon\lambda/3}.
  \end{displaymath}\end{it}

{\sc Proof:}
By the Cauchy-Schwarz inequality, we have

  \begin{equation}
  S_{\theta}(u,\sigma)\ll
  {\bf N}(T)^{1/2}\left(\hspace{-0.4cm} 
  \sum\limits_{\scriptsize \begin{array}{cccc} &\rho:\ \! 0< \gamma \le T,&\\
  &\sigma\le \beta\le \sigma+1/\log N& \end{array}} \hspace{-0.4cm}
  \left\vert D_{\theta}(u,i\gamma)\right\vert^2\right)^{1/2},\label{29}
  \end{equation}
where ${\bf N}(T)$ denotes the number of all non-trivial zeta zeros $\rho$
with $0<\gamma\le T$. By the
well-known properties of the set of zeta zeros, we can split the set
of zeros $\rho$ satisfying the conditions $0< \gamma \le T,$
$\sigma\le \beta\le \sigma+1/\log N$ into $O(\log T)$ subsets ${\bf S}$ 
satisfying the condition

  \begin{displaymath}
  \rho_1, \rho_2\in {\bf S},\ \rho_1\not=\rho_2 \ \Longrightarrow \
  \vert \mbox{ Im } \rho_1 - \mbox{ Im } \rho_2 \ \! \vert \ge 1.
  \end{displaymath}
Employing {\sc Proposition 2}, we get 

  \begin{equation} 
  \sum\limits_{\rho\in {\bf S}}
  \left\vert D_{\theta}(u,i\gamma)\right\vert^2 \ll 
  (T+N^{\lambda})N^{-\lambda}(\log N)\label{30}
  \end{equation}
for $N^{\lambda}<u\le (2N)^{\lambda}$.

Combining (\ref{107}), (\ref{29}), (\ref{30}) 
and ${\bf N}(T)\ll T\log T$, and taking the condition $\varepsilon \le 
f(\lambda,h)/\lambda$ of {\sc Lemma 2} into account, we obtain the
desired bound. $\Box$\\
 
{\sc Lemma 3:} \begin{it} For
$1-\Delta<\sigma\le 1-1/(\log N)$ we have   

  \begin{displaymath}
  N^{\sigma}S_{\theta}(u,\sigma)\ll
  N\exp\left(-(\log N)^{1/4}\right).
  \end{displaymath}
\end{it}

{\sc Proof:} 
By {\sc Proposition 3}, there is no zeta zero $\rho$ with $0<\gamma\le T$ 
on the right side of the line Re $s=\kappa(T)$, where

  \begin{displaymath}
  \kappa(T):= 1-
  C(\log T)^{-2/3}(\log \log T)^{-1/3}.
  \end{displaymath} 
Therefore, we can assume that $\sigma\le \kappa(T)$.

We first consider the case when $h\ge 2$. 
From the trivial estimate

  \begin{displaymath}
  D_{\theta}(u,i\gamma)\ll 1,
  \end{displaymath}
the second zero density estimate of
{\sc Proposition 6} and (\ref{100}),
we obtain

  \begin{equation}
  N^\sigma S_{\theta}(u,\sigma)\ll
  N^{\sigma+(15(1-\sigma))/(11(3\sigma-1))}(\log N)^{44}.\label{31}
  \end{equation}
We notice that 

  \begin{equation}
  \frac{15}{11(3\sigma-1)}< 1-\frac{1}{3751}\label{32}
  \end{equation}
if $\sigma> 1-\Delta=0.788$. From (\ref{31}), (\ref{32}) and the above 
assumption $\sigma\le \kappa(T)$ follows

  \begin{displaymath} 
  N^{\sigma} S_{\theta}(u,\sigma) \ll N^{1-(1-\kappa(T))/3751}(\log N)^{44}.
  \end{displaymath}
From this, we obtain

  \begin{equation} 
  N^{\sigma} S_{\theta}(u,\sigma) \ll N\exp\left(-(\log N)^{1/4}\right)
  \label{33}
  \end{equation}
by a short 
calculation. This completes the proof for the case when $h\ge 2$. 

Now, let $h=1$. Then, similar to (\ref{31}), we get

  \begin{equation}
  N^\sigma S_{\theta}(u,\sigma)\ll
  N^{\sigma+(15(1-\sigma))/(8(3\sigma-1))}(\log N)^{44}\label{34}
  \end{equation}
by using (\ref{106}).
We notice that 

   \begin{equation}
  \frac{15}{8(3\sigma-1)}< 1-\frac{1}{46} \label{35}
  \end{equation}
if $\sigma>1-\Delta=35/36$. In a similar manner like in the case when $h\ge 2$,
from (\ref{34}) and (\ref{35}), we anew obtain (\ref{33}). This 
completes the proof.      
$\Box$\\

\section{Zero-detection method with weights}
Next, we use a modified form of the zero detection method to handle
the sum $S_{\theta}(u,\sigma)$ in the range $6/11<\sigma\le 1-\Delta$.

We note that by (\ref{18}) the condition $\log N \ll \log T \ll \log N$ of
{\sc Proposition 4} is satisfied if $\varepsilon<1$.\\

{\sc Lemma 4:} \begin{it} Suppose that $1/2\le \sigma \le 1-1/(\log N)$.
Then, under the conditions and using the definitions of
{\sc Proposition} {\rm 4}, we have

  \begin{displaymath}
  S_{\theta}(u,\sigma)\ \ll \ \left(V_1(u,\sigma)^{1/(2h)}+
  V_2(u,\sigma)^{1/2}\right)N^{\varepsilon\lambda/400}+(\log N)^3
  \end{displaymath}
with
  \begin{eqnarray}\hspace*{0.7cm}
  V_1(u,\sigma)&:=& {\bf N}(\sigma,T)^{2h-3/2} T^{1/2} 
  Y^{1-2\sigma}\cdot  \label{36} \\ & &
  \int\limits_{-(\log N)^2}^{(\log N)^2}
  \sum\limits_{\rho} {}^{(\sigma)} 
  \left\vert M_X(1/2+i(\gamma+t))\right\vert^2 \cdot 
  \left\vert D_{\theta}(u,i\gamma)\right\vert^{2h} 
  \ {\rm d}t \nonumber
  \end{eqnarray}\\
and 
  \begin{eqnarray}
  & & \label{37} \\ & & V_2(u,\sigma) \nonumber\\  
  &:=& {\bf N}(\sigma,T) 
  \sup\limits_{X<v\le Y(\log N)^2}
  \sum\limits_{\rho} {}^{(\sigma)}
  \left\vert \sum\limits_{X<k\le v} 
  a(k)e^{-k/Y} k^{-(\sigma+i\gamma)} \right\vert^2 \cdot
  \left\vert D_{\theta}(u,i\gamma) \right\vert^2,\nonumber 
  \end{eqnarray}
where the notation $(\sigma)$ attached to the summation symbol indicates the
summation condition ``\ $(\log N)^2<\gamma\le T$ and $\sigma\le \beta\le
\sigma+1/(\log N)$''. The term $V_2(u,\sigma)$ is supposed to equal 0 if
$X\ge Y(\log N)^2$.\end{it}\\

{\sc Proof:} We first consider the contribution of zeta zeros with
small ima\-ginary part $\gamma\le (\log N)^2$. By
${\bf N}((\log N)^2)\ll (\log N)^3$ and the trivial estimate
$D_{\theta}(u,i\gamma)\ll 1$, we get

  \begin{displaymath}
  \hspace{-0.4cm} \sum\limits_{\scriptsize \begin{array}{cccc} &\rho:\ \!
  0< \gamma \le (\log N)^2,&\\
  &\sigma\le \beta\le\sigma+1/\log N& \end{array}} \hspace{-0.4cm}
  \vert D_{\theta}(u,i\gamma) \vert \ll (\log N)^3.
  \end{displaymath}   
It remains to prove that

  \begin{equation}
  \sum\limits_{\rho} {}^{(\sigma)} \ \vert D_{\theta}(u,i\gamma) \vert \ll 
  \left(V_1(u,\sigma)^{1/(2h)}+
  V_2(u,\sigma)^{1/2}\right)N^{\varepsilon\lambda/400}. \label{38}
  \end{equation}
By {\sc Proposition 4}, for every sufficiently large $N$ we have

  \begin{equation}
  \sum\limits_{\rho} {}^{(\sigma)} \ \vert D_{\theta}(u,i\gamma) \vert \le
  \left(\sum\limits_{\rho} {}^{(\sigma,1)}+ 
  \sum\limits_{\rho} {}^{(\sigma,2)}\right)\
  \left\vert D_{\theta}(u,i\gamma)\right\vert, \label{39}
  \end{equation}
where the notation $(\sigma,j)$ attached to the summation 
symbol on the right side
 indicates that
$\rho$ is a zero of type $j$ (for the definition of ``type $j$''
see {\sc Proposition 4}) satisfying the summation condition 
$(\sigma)$, {\it i.e.}
$(\log N)^2<\gamma\le T$ and $\sigma\le \beta\le \sigma+1/(\log N)$.
We now separately consider the two sums on the right side of (\ref{39}).

By H\" older`s inequality and 
$1/2<U_1(\rho)$ for every zero $\rho$ of type $1$, we get
  \begin{equation}
  \left(\sum\limits_{\rho} {}^{(\sigma,1)}\
  \left\vert D_{\theta}(u,i\gamma)\right\vert \right)^{h}\ll
  {\bf N}(\sigma,T)^{h-1}
  \sum\limits_{\rho} {}^{(\sigma)}\ U_1(\rho)\cdot  
  \left\vert D_{\theta}(u,i\gamma)\right\vert^{h}.\label{40}
  \end{equation}
Using the triangle and the Cauchy-Schwarz inequality, and applying
Stirling$^{\prime}$s formula to the Gamma factor contained in the integrand
on the right side of (\ref{13}), we derive
  \begin{eqnarray}
  & & \label{41} \\ & & \left( \sum\limits_{\rho} {}^{(\sigma)}\ U_1(\rho)  
  \vert D_{\theta}(u,i\gamma) \vert^{h}\right)^2 \nonumber\\ &\ll&
  (\log N)^2 
  {\bf N}(\sigma,T)^{1/2}Y^{1-2\sigma} \left(\sum\limits_{\rho} {}^{(\sigma)}
  \int\limits_{-(\log N)^2}^{(\log N)^2} \vert \zeta(1/2+i(\gamma+t))
  \vert^4\ {\rm d}t\right)^{1/2} \cdot \nonumber\\
  & & \left(\int\limits_{-(\log N)^2}^{(\log
  N)^2} \sum\limits_{\rho} {}^{(\sigma)} \ \vert 
  M_X(1/2+i(\gamma+t))\vert^2 \cdot 
  \vert D_{\theta}(u,i\gamma)\vert^{2h}\ {\rm d}t\right). \nonumber 
  \end{eqnarray}
Taking notice of $\log T\ll \log N$ and 
${\bf N}(t+1)-{\bf N}(t)=O(\log(2t))$ for $t\ge 1$, and using
{\sc Proposition 5}, the term involving the
$\zeta$-function on the right side
can be estimated by
  \begin{eqnarray}
  & & \sum\limits_{\rho} {}^{(\sigma)}\
  \int\limits_{-(\log N)^2}^{(\log N)^2} \vert \zeta(1/2+i(\gamma+t))
  \vert^4\ {\rm d}t\label{42}\\ &\ll& (\log N)^3 
  \int\limits_{0}^{T+(\log N)^2} \vert \zeta(1/2+it) \vert^4\ {\rm d}t
  \nonumber\\ &\ll& T(\log N)^7.\nonumber
  \end{eqnarray}

Since    
$1/2< U_2(\rho)$ for every zero $\rho$ of type $2$, we have

  \begin{displaymath}
  \sum\limits_{\rho} {}^{(\sigma,2)}\
  \left\vert D_{\theta}(u,i\gamma)\right\vert 
  \ll \sum\limits_{\rho} {}^{(\sigma)}\ U_2(\rho) 
  \left\vert D_{\theta}(u,i\gamma)\right\vert.
  \end{displaymath}
From that, applying partial summation to the term $U_2(\sigma)$
on the right side and 
taking $N^{1/(\log N)}=e$ into account, we obtain
  \begin{eqnarray}
  & & \sum\limits_{\rho} {}^{(\sigma,2)}\
  \left\vert D_{\theta}(u,i\gamma)\right\vert \label{43}\\
  &\ll& \sup\limits_{X<v\le Y(\log N)^2}
  \ \sum\limits_{\rho} {}^{(\sigma)}\
  \left\vert \sum\limits_{X<k\le v} 
  a(k)e^{-k/Y} k^{-(\sigma+i\gamma)} \right\vert \cdot 
  \left\vert D_{\theta}(u,i\gamma)\right\vert.\nonumber
  \end{eqnarray}
By the Cauchy-Schwarz inequality, we get
 \begin{eqnarray}
  & & \left(\sum\limits_{\rho} {}^{(\sigma)}\
  \left\vert \sum\limits_{X<k\le v} 
  a(k)e^{-k/Y} k^{-(\sigma+i\gamma)} \right\vert \cdot 
  \left\vert D_{\theta}(u,i\gamma)\right\vert\right)^2 \label{44} \\ &\ll&
  {\bf N}(\sigma,T)\sum\limits_{\rho} {}^{(\sigma)}\
  \left\vert \sum\limits_{X<k\le v} 
  a(k)e^{-k/Y} k^{-(\sigma+i\gamma)} \right\vert^2 \cdot 
  \left\vert D_{\theta}(u,i\gamma)\right\vert^{2}.\nonumber
  \end{eqnarray}

Combining (\ref{39})-(\ref{44}), we obtain (\ref{38}). 
This completes the proof.$\Box$\\  

To bound $V_j(u,\sigma)$ $(j=1,2)$ by 
simple terms involving the parameters $X$ and $Y$,  
we apply {\sc Theorem 5} after splitting the set of zeta zeros satisfying 
the condition $(\sigma)$
into $O(\log T)$ subsets ${\bf S}$ such that
 $\vert\!$ Im $\rho_1-$ Im
$\rho_2\vert\ge 1$ for every pair $\rho_1,\rho_2\in {\bf S}$ with 
$\rho_1\not=\rho_2$. We take into consideration that 

  \begin{displaymath}
  \vert a(k) \vert \le \tau(k) \ll N^{\varepsilon_1}
  \end{displaymath} 
for $X<k\le Y(\log N)^2$, 
where $\tau(k)$ denotes the number of divisors of $k$ and 
$\varepsilon_1$ is any positive constant. Moreover,
we point out that for irrational $\theta$ the additional conditions 
$u=L_2\le T^{1/2}$ and $h=1$ in {\sc Theorem 5} are really 
satisfied. Indeed, at the beginning of section 3 we noticed that $h=1$ if 
$\theta$ is irrational, and 
by (\ref{109}) and  the condition $\lambda< 5/19$ in 
{\sc Theorem 3}, the inequality 
$u\le T^{1/2}$ is satisfied if 
$\varepsilon \le (5-19\lambda)/(12\lambda)$ and $N\ge 5$ 
(the two latter conditions
may be supposed without loss of generality).

In this manner, we obtain the following result.\\

{\sc Lemma 5:} \begin{it} Suppose that $1/2\le \sigma\le 1-1/\log N$.
If $\theta$ is irrational, 
then suppose that $\lambda<5/19$ and 
$\varepsilon \le (5-19\lambda)/(12\lambda)$. Then,
on the conditions of {\sc Proposition} {\rm 4}, we have

  \begin{displaymath}
  V_1(u,\sigma)\ \ll\ {\bf N}(\sigma,T)^{2h-3/2} T^{1/2}
  Y^{1-2\sigma}\left(TN^{-h\lambda}+X\right) N^{h\varepsilon\lambda/200}
  \end{displaymath} 
and

  \begin{displaymath}
  V_2(u,\sigma) \ \ll\ {\bf N}(\sigma,T)
  \left(TN^{-\lambda}X^{1-2\sigma}+Y^{2(1-\sigma)}\right)
  N^{\varepsilon\lambda/200},
  \end{displaymath}
the implied $\ll$-constant only depending on $\varepsilon$.\end{it}
 
\section{Estimation of $N^{\sigma}S_{\theta}(u,\sigma)$ for
  $6/11<\sigma\le 1-\Delta$}

The final step of the proof of {\sc Theorem 3} is to show\\

{\sc Lemma 6:} \begin{it} The estimate $(\ref{28})$ holds true for
  $6/11<\sigma\le 1-\Delta$. \end{it} \\

{\sc Lemma} 6 follows from the preceding {\sc Lemmas} 4, 5 and the following\\

{\sc Lemma 7:} \begin{it} Without loss of generality assume that 
$\varepsilon\le 1/(10\lambda)$.
Then for any $\sigma$ in the range 
$6/11< \sigma\le 1-\Delta$ there are parameters $X$, $Y$ satisfying the 
conditions of {\sc Proposition} {\rm 4}, such that
  
  \begin{equation}
  R_j(\sigma)\ll N^{2-\varepsilon\lambda/90} \ \ \ \ \ \ (j=1,...,4)\label{45}
  \end{equation}
as $N\rightarrow \infty$, where   
  \begin{eqnarray*}
  R_1(\sigma)&:=&{\bf N}(\sigma,T)^{2-3/(2h)} 
  T^{3/(2h)}N^{2\sigma-\lambda}Y^{(1-2\sigma)/h},
  \nonumber\\ \nonumber\\
  R_2(\sigma)&:=&{\bf N}(\sigma,T)^{2-3/(2h)} 
  T^{1/(2h)}N^{2\sigma}X^{1/h}Y^{(1-2\sigma)/h}, \nonumber\\ \nonumber\\
  R_3(\sigma)&:=&{\bf
  N}(\sigma,T)TN^{2\sigma-\lambda}X^{1-2\sigma},\nonumber\\ \nonumber\\
  R_4(\sigma)&:=&{\bf N}(\sigma,T)N^{2\sigma} Y^{2(1-\sigma)},
  \end{eqnarray*}
the implied $\ll$-constant in $(\ref{45})$ depending only on $\varepsilon$.
\end{it}\\

{\sc Proof:} Firstly, we consider the case when $h=1$.
We put

  \begin{displaymath}
  Y:=N^{1-\varepsilon\lambda/5}(1+{\bf N}(\sigma,T))^{-1/(2(1-\sigma))}
  \end{displaymath}
and 
 
  \begin{displaymath}
  X:=1+Y^{2\sigma-1}N^{2(1-\sigma)(1-2\varepsilon\lambda/5)}T^{-1/2}
  (1+{\bf N}(\sigma,T))^{-1/2}.
  \end{displaymath}

We now derive some simple estimates for $X$ and $Y$ to verify the
conditions of {\sc Proposition 4}. Trivially, we have $1\le X$. 
By (\ref{106}) and the well-known bound

  \begin{displaymath}
  {\bf N}(\sigma,T)\ll T^{12(1-\sigma)/5}(\log T)^{44} 
  \end{displaymath} 
following from {\sc Proposition 6}, we get
    
  \begin{displaymath} 
  N^{(1-\lambda)/2} \ll Y \ll N.
  \end{displaymath}
This implies $\log N \ll \log Y \ll \log N$. Consequently, $\log T\ll
\log Y \ll \log T$. The last condition to be verified is $\log X\ll \log N$. 
To prove this inequality, it suffices to show that 
$X\le Y$ for sufficiently large $N$, which latter   
follows from $1=o(Y)$ and 

  \begin{equation}     
  Y^{2\sigma-1}N^{2(1-\sigma)(1-2\varepsilon\lambda/5)}T^{-1/2}
  (1+{\bf N}(\sigma,T))^{-1/2}=o(Y)\label{47}
  \end{equation}
as $N\rightarrow\infty$. The bound (\ref{47}) can be easily obtained from the
definition of $Y$, the bound ${\bf N}(\sigma,T)\ll T(\log T)$ and 
$\sigma\le 1-\Delta=35/36$.
Therefore, all conditions of {\sc Proposition 4} to $X$, $Y$ are satisfied.
 
Next, we calculate the order of magnitude of the terms $R_j(\sigma)$. 
From the definitions of $X$, $Y$ and $\sigma\le 35/36$, 
we obtain

  \begin{equation}
  R_4(\sigma)\ll\ N^{2-\varepsilon\lambda/90},
  \label{48}
  \end{equation}

  \begin{equation}
  R_2(\sigma)=N^{2-\varepsilon\lambda/45}+
  R_1(\sigma)T^{-1}N^{\lambda}
  \label{49} 
  \end{equation}
 and

  \begin{equation}
  R_1(\sigma)\ll (1+{\bf N}(\sigma,T))^{\sigma/(2(1-\sigma))}
  T^{3/2}N^{1-\lambda+\varepsilon\lambda/5}.
  \label{50}
  \end{equation}
From {\sc Proposition 6}, we derive
  \begin{eqnarray}
  & & \sup\limits_{1/2\le \sigma \le 35/36}
  (1+{\bf N}(\sigma,T))^{\sigma/(1-\sigma)}\label{51}\\ 
  &\ll& T^{\varepsilon\lambda/5}\left(\sup\limits_{1/2\le \sigma \le 3/4} 
  T^{3\sigma/(2-\sigma)} +
  \sup\limits_{3/4\le \sigma \le 35/36}  T^{3\sigma/(3\sigma-1)}
  \right).\nonumber
  \end{eqnarray}
The
function $g_1(\sigma):=\sigma/(2-\sigma)$ is monotone increasing on the
interval $[1/2,3/4]$, and the function
$g_2(\sigma):=\sigma/(3\sigma-1)$ is
monotone decreasing on the interval $[3/4,35/36]$. Therefore, from 
(\ref{51}) follows

  \begin{equation}
  \sup\limits_{1/2\le \sigma \le 35/36}
  (1+{\bf N}(\sigma,T))^{\sigma/(1-\sigma)} \ll T^{9/5+\varepsilon\lambda/5}.
  \label{52}
  \end{equation}
Combining the first line of 
(\ref{106}), (\ref{50}) and (\ref{52}), we get 

  \begin{equation} 
  R_1(\sigma)\ll N^{2-\varepsilon\lambda/2}.\label{53}
  \end{equation}
From $N^{\lambda}\le T$, (\ref{49}) and (\ref{53}), we obtain

   \begin{equation}
   R_2(\sigma)\ll
   N^{2-\varepsilon\lambda/45}.\label{54}
   \end{equation}

The last step is to verify the bound

  \begin{equation}
  R_3(\sigma)\ll N^{2-\varepsilon\lambda/5}.
  \label{55}
  \end{equation}
From $X\le Y$ (which we have seen above) and the definitions of $X$ and $Y$, 
we conclude

  \begin{displaymath}
  \left(\frac{Y}{X}\right)^{2\sigma-1}\le \frac{Y}{X} \le 
  \left(\frac{T}{1+{\bf N}(\sigma,T)}\right)^{1/2}\cdot 
  N^{\varepsilon\lambda/5},
  \end{displaymath}
from which follows 

  \begin{eqnarray*}
  R_3(\sigma)&=& 
  {\bf N}(\sigma,T) TN^{2\sigma-\lambda}X^{1-2\sigma}\\ &\le&
  {\bf N}(\sigma,T)^{1/2} T^{3/2} N^{2\sigma-\lambda}Y^{1-2\sigma}
  N^{\varepsilon\lambda/5}\\ &=& 
  R_1(\sigma)N^{\varepsilon\lambda/5}.
  \end{eqnarray*}
Combining this inequality and (\ref{53}), we get (\ref{55}). 

By (\ref{48}), (\ref{53}), (\ref{54}) and (\ref{55}), 
the bound (\ref{45}) is satisfied for $j=1,...,4$. This 
completes the proof for the case when $h=1$.

Secondly, we consider the case when $h\ge 2$. We observe that 
$N^{h\lambda}\le T$ by 
(\ref{9}), (\ref{101}) and the assumption 
$\varepsilon\le 1/(10\lambda)$ of {\sc Lemma 7}.
Here we put
  
  \begin{displaymath}
  X:=TN^{-h\lambda}
  \end{displaymath}
and

  \begin{displaymath}
  Y:=N^{1-\varepsilon\lambda/4}(1+{\bf N}(\sigma,T))^{-1/(2(1-\sigma))}
  \end{displaymath}
unlike in the case when $h=1$.
Using {\sc Proposition 6}, (\ref{100}) and $\varepsilon\le 1/(10\lambda)$,
it is easily verified that $X$ and $Y$ satsify the conditions of 
{\sc Proposition 4}. Further, it is an immediate consequence of the 
definition of $Y$ and the condition $6/11<\sigma\le 1-\Delta=
0.788$ that
(\ref{45}) holds true for $j=4$. Here, as in the following, we 
use the condition $6/11<\sigma\le
0.788$ in order to obtain the correct $\varepsilon$-terms.

From the definitions of $X$ and $Y$ follows
  \begin{eqnarray}
  & & R_1(\sigma)+R_2(\sigma)\label{56} \\ &\ll&  
  {\bf N}(\sigma,T)^{2-3/(2h)+(2\sigma-1)/(2h(1-\sigma))}
  T^{3/(2h)}N^{2\sigma(1-1/h)+1/h-\lambda+\varepsilon\lambda/(7h)}.\nonumber
  \end{eqnarray}
By {\sc Proposition 6}, we have

  \begin{equation}
  {\bf N}(\sigma,T)\ll T^{A(\sigma)(1-\sigma)+\varepsilon\lambda/10},\label{57}
  \end{equation}
where 

  \begin{displaymath}
  A(\sigma)=\left\{\hspace{-0.4cm} 
  \begin{array}{llll} &3/(2-\sigma)\ \mbox{ if }\
  1/2\le \sigma\le 3/4,& \\ \\ &3/(3\sigma-1)\ \mbox{ if }\ 3/4<
  \sigma\le 1.& \end{array} \right.
  \end{displaymath}
Combining (\ref{18}), (\ref{56}) and (\ref{57}), and taking (\ref{10}) and 
$h\ge 2$ 
into consideration, we get

  \begin{equation}
  R_1(\sigma)+R_2(\sigma)\ll N^{r_1(\sigma)+c_1-\varepsilon\lambda/10},
  \label{58}
  \end{equation}
where
  \begin{displaymath}
  r_1(\sigma):=\left(\lambda+f(\lambda,h)\right)
  \left(2-\frac{2}{h}+
  \left(\frac{5}{2h}-2\right)\sigma\right)A(\sigma)+
  2\left(1-\frac{1}{h}\right)\sigma
  \end{displaymath}
and 

  \begin{displaymath}
  c_1:=\left(\lambda+f(\lambda,h)\right)\cdot 
  \frac{3}{2h}+\frac{1}{h}-\lambda.
  \end{displaymath}

Our next aim is to show that $r_1(\sigma)$ is monotone increasing
on the interval $6/11<\sigma< 3/4$ and monotone decreasing on the
interval  $3/4< \sigma\le 0.788$. For $6/11< \sigma<3/4$ we have

  \begin{displaymath}
  r_1^{\prime}(\sigma)=-
  \left(\lambda+f(\lambda,h)\right)
  \left(2-\frac{3}{h}\right)\cdot \frac{3}{(2-\sigma)^2}+
  2\left(1-\frac{1}{h}\right).
  \end{displaymath}
From that, (\ref{10}) and $h\ge 2$, we obtain

  \begin{displaymath}
  r_1^{\prime}(\sigma) \ge \frac{14}{55}+\frac{34}{55h}> 0
  \end{displaymath}
for $6/11< \sigma<3/4$. Hence, $r_1(\sigma)$ is monotone increasing
on this interval.
For $3/4< \sigma\le 0.788$ we have

  \begin{displaymath}
  r_1^{\prime}(\sigma)=
  -\left(\lambda+f(\lambda,h)\right)
  \left(12-\frac{21}{2h}\right)\cdot \frac{1}{(3\sigma-1)^2}+
  2\left(1-\frac{1}{h}\right).
  \end{displaymath}
From that and (\ref{10}), we obtain

  \begin{displaymath}
  r_1^{\prime}(\sigma)<-0.6+\frac{0.4}{h}< 0
  \end{displaymath}
for $3/4< \sigma\le 0.788$. Hence, $r_1(\sigma)$ is monotone decreasing
on this interval.

We note that the function $r_1(\sigma)$ is continuous on the
interval $(6/11,0.788]$ since $A(\sigma)$ is continuous on this
interval. From that and the above observations,
we conclude that the exponent $r_1(\sigma)+c_1-\varepsilon\lambda/10$
on the right side of (\ref{58}) takes
its maximum at the point $\sigma_0=3/4$. Furthermore, from  

  \begin{displaymath}
  \lambda+f(\lambda,h)\le \lambda+f_1(\lambda,h)= 
  \frac{5}{12}+\frac{5h\lambda}{6(h+1)},
  \end{displaymath}
we obtain

  \begin{displaymath}
  r_1\left(3/4\right)+c_1\le 2
  \end{displaymath}
by a short calculation. From that and (\ref{58}), 
we derive (\ref{45}) for $j=1,2$.

Finally, we evaluate the term $R_3(\sigma)$.
From (\ref{18}), (\ref{57}) and the definition of $X$, we obtain
  
  \begin{equation}
  R_3(\sigma)\ll N^{r_2(\sigma)-(h+1)\lambda-\varepsilon\lambda/2},\label{59}
  \end{equation}
where 

  \begin{displaymath}
  r_2(\sigma):=
  \left(\lambda+f(\lambda,h)\right)(2+A(\sigma))
   (1-\sigma)+2(1+h\lambda)\sigma.
  \end{displaymath}
For $6/11< \sigma<3/4$ we have

  \begin{displaymath}
  r_2^{\prime}(\sigma)=-
  \left(\lambda+f(\lambda,h)\right)
  \left(2+\frac{3}{(2-\sigma)^2}\right)+
  2\left(1+h\lambda\right).
  \end{displaymath}
From that and (\ref{10}), we obtain

  \begin{displaymath}
  r_2^{\prime}(\sigma) > 0
  \end{displaymath}
for $6/11< \sigma<3/4$. Hence, $r_2(\sigma)$ is monotone increasing
on this interval. At the end of this section, 
we shall separately prove that $r_2(\sigma)$ is monotone decreasing on the 
interval
$3/4< \sigma\le 0.788$. 

Like $r_1(\sigma)$, the function $r_2(\sigma)$ is continuous on the
interval $(6/11,0.788]$. Consequently, the exponent 
$r_2(\sigma)-(h+1)\lambda-\varepsilon\lambda/2$
on the right side of (\ref{59}) takes
its maximum at the point $\sigma_0=3/4$. Furthermore, from  

  \begin{displaymath}
  \lambda+f(\lambda,h)\le \lambda+f_2(\lambda,h)=
  \frac{5}{11}+\frac{(10-5h)\lambda}{11},
  \end{displaymath}
we obtain

  \begin{displaymath}
  r_2\left(3/4\right)-(h+1)\lambda\le 2
  \end{displaymath}
by a short calculation. From that and (\ref{59}), 
we derive (\ref{45}) for $j=3$.
This completes the proof of {\sc Lemma 7}.
$\Box$\\

By proving {\sc Lemma 7} we have also completed the proof of {\sc Theorem
3}. 

It remains to show that $r_2^{\prime}(\sigma)<0$  for
$3/4< \sigma\le 0.788$ if $h\ge 2$. On this interval, we have

  \begin{displaymath}
  r_2^{\prime}(\sigma)=
  -\left(\lambda+f(\lambda,h)\right)
  \left(2+\frac{6}{(3\sigma-1)^2}\right)+
  2\left(1+h\lambda\right).
  \end{displaymath}
Thus, $r_2^{\prime}(\sigma)<0$ for
$3/4< \sigma\le 0.788$ is equivalent to 

  \begin{displaymath}
  2\left(2+\frac{6}{1.364^2}\right)^{-1}=:\xi<
  \frac{\lambda+f(\lambda,h)}{1+h\lambda},
  \end{displaymath}
where we have $\xi\approx 0.3828$.

By definition, for $k\in {\bf N}$ we have

  \begin{displaymath}
  \frac{\lambda+f(\lambda,k)}{1+k\lambda}=
  \min\left\{\frac{\lambda+f_{1}(\lambda,k)}{1+k\lambda}, 
  \frac{\lambda+f_{2}(\lambda,k)}{1+k\lambda}\right\} 
  \end{displaymath}
with 
  
  \begin{displaymath}
  \frac{\lambda+f_{1}(\lambda,k)}{1+k\lambda}=\frac{5}{12}\cdot 
  \left(1-\left(1-\frac{1}{1+k\lambda}\right)\left(1-\frac{2}{1+k}
  \right)\right) 
  \end{displaymath}
and 
  
  \begin{displaymath}
  \frac{\lambda+f_{2}(\lambda,k)}{1+k\lambda}=\frac{5}{11}\cdot 
  \frac{1+(2-k)\lambda}{1+k\lambda}. 
  \end{displaymath}
For fixed $k\ge 2$ the functions 
$g_{i}(\lambda):=(\lambda+f_{i}(\lambda,k))/(1+k\lambda)$ 
$(i=1,2)$ are obviously monotone decreasing for $\lambda>0$. Furthermore, 
by (\ref{8}) and $h\ge 2$, we have

  $$
  \lambda\le \phi_{1,h}=\frac{1}{2(6h-1-11/h)}
  $$
and 

  $$
  f(\phi_{1,h},h)=f_2(\phi_{1,h},h).
  $$
Hence, it suffices to prove that

  \begin{equation}
  \xi<z_h:=\frac{5}{11}\cdot 
  \frac{1+(2-h)/(2(6h-1-11/h))}{1+h/(2(6h-1-11/h))}. \label{200}
  \end{equation}
The inequality (\ref{200}) holds true for $h=2,3,4$. Furthermore, 
it is easily seen that the sequence $\left(z_h\right)$ is monotone
decreasing for $h\ge 4$, and we have 

 $$
 \lim\limits_{h\rightarrow\infty} z_h=\frac{5}{13}>\xi.
 $$
This completes the proof. $\Box$
  
\section{Proofs of Theorems 4, 4$^*$ and 5}
{\sc Theorem} $4^*$ can be derived from {\sc Theorem 4} in a standard way 
using the inequality

  \begin{displaymath}
  \vert f(x) \vert \le \int\limits_{x-1/2}^{x+1/2} (\vert f(t) \vert + 
  \vert f^{\prime}(t)\vert )\ {\rm d}t
  \end{displaymath}
which is valid for every continuously differentable function $f: [x-1/2,x+1/2]
\rightarrow {\bf C}$.\\

To derive {\sc Theorem 5} from {\sc Theorem} 
$4^*$, we proceed as follows: In case $\theta$
is rational we write the shifted Dirichlet polynomial on the left side of
$(\ref{17})$ as an ordinary one via the relation 

$$
\sum\limits_{L_1<l\le L_2} b_{l}(l+\theta)^{i\alpha t_r-1}
=q^{-i\alpha t_r+1} \sum\limits_{L_1<l\le L_2}
b_l(ql+m)^{i\alpha t_r-1},
$$ 
where $\theta=m/q$, $m$ and $q$ being non-negative integers. We then write the 
$2h$-th power of the absolute value of 
the Dirichlet polynomial on the right side in the form

  \begin{displaymath}
  \left\vert \sum\limits_{L_1<l\le L_2}
  b_l(ql+m)^{i\alpha t_r-1} \right\vert^{2h} =
  \left\vert 
  \sum\limits_{(qL_1+m)^h<n\le (qL_2+m)^h}
  c_n n^{i\alpha t_r} \right\vert^2,
  \end{displaymath}
where 

  \begin{displaymath}
  c_n:=n^{-1}\hspace{-0.4cm} 
  \sum\limits_{\scriptsize \begin{array}{cccc}&L_1<l_1,...,l_h\le L_2,& \\
  &n=(ql_1+m)\cdots (ql_h+m)& 
  \end{array}}  \hspace{-0.4cm} 
  b_{l_1}\cdots b_{l_h}.
  \end{displaymath} 
We note that
  
  $$
  \vert c_n\vert \le B^hn^{-1}d_h(n) \ll B^hn^{\varepsilon_2-1},
  $$
where $d_h(n)$ denotes
the divisor function of order $h$ and $\varepsilon_2$ is any
positive constant. Now, we divide each of the Dirichlet polynomials

$$
\sum\limits_{K_1<k\le K_2} 
  a_{k}k^{-(\sigma+it_r)}
$$
and 

$$
\sum\limits_{(qL_1+m)^h<n\le (qL_2+m)^h} c_n n^{i\alpha t_r}
$$
into $O(\log K_2)$ and $O(\log L_2)$ 
partial sums over ranges of the form $K<k\le 2K$ and $L<l\le 2L$ 
respectively, use the Cauchy-Schwarz inequality, multiply out the two 
resulting 
sums of squares of absolute values of Dirichlet polynomials, and sum up over
$r$. In this manner, we obtain 
a sum of terms having the same shape as the one on the left side of
(\ref{16}), where now $\theta=0$. Applying {\sc Theorem} $4^*$ with 
$\theta=0$ to these terms, we obtain the desired bound.

When $\theta$ in {\sc Theorem 5} is irrational, it is supposed that $h=1$.
Now, we just split up the ordinary and the shifted Dirichlet polynomial 
on the left side of (\ref{17}) in the same manner as above, use the 
Cauchy-Schwarz inequa\-lity, multiply out, sum up over $r$ and 
apply {\sc Theorem} $4^*$. In this way, we again obtain the desired bound.
This completes the
proof of {\sc Theorem 5}.\\ 
      
We now turn to proving {\sc Theorem} 4. If $\theta=0$, {\sc Theorem 4} is 
nothing but a slight modification of Theorem 1 in [BaH]. However,
in the case when $\theta\not= 0$ {\sc Theorem 4} actually appears to be a new 
result, which we shall prove in the following. 

Without loss of
generality, we assume that $A=B=1$ and $\alpha>0$. 
We denote the integral in question
on the left side of (\ref{15}) by $I$.

Multiplying out the integrand contained in $I$, 
integrating the resulting
fourfold sum term by term and using the standard inequalities

  \begin{displaymath}
  \int\limits_{0}^{T} x^{it}\ {\rm d}t \ll \min\{T,\vert \log x \vert^{-1}\}
  \end{displaymath}
for $x>0$ and 

  \begin{displaymath}
  \vert \log \omega \vert \ \gg \ \vert\omega-1\vert 
  \end{displaymath}
for $2^{-(\alpha+1)}\le \omega\le 2^{\alpha+1}$, we obtain
  \begin{eqnarray} 
  I &\le& \sum\limits_{k_1,k_2\sim K} \ \sum\limits_{l_1,l_2\sim L}
  \min\left\{T,\ \left\vert 
  \frac{k_1(l_2+\theta)^{\alpha}}{k_2(l_1+\theta)^{\alpha}}-1
  \right\vert^{-1} \right\} \label{60} \\ &\le& 2^{\alpha} 
  \sum\limits_{k_1,k_2\sim K} \ \sum\limits_{l_1,l_2\sim L}
  \min\left\{T,\ \left\vert \frac{k_1}{k_2}-
  \left(\frac{l_1+\theta}{l_2+\theta}\right)^{\alpha} \right\vert^{-1}\right\}
  \nonumber\\ &\le& 2^{\alpha} \sum\limits_{d\le 2K} \hspace{-0.4cm}
  \sum\limits_{\scriptsize \begin{array}{cccc} &k_1,k_2\sim K/d,&\\ 
  &(k_1,k_2)=1& \end{array}} \hspace{-0.4cm} \sum\limits_{l_1,l_2\sim L}
  \min\left\{T,\ \left\vert \frac{k_1}{k_2}-
  \left(\frac{l_1+\theta}{l_2+\theta}\right)^{\alpha} \right\vert^{-1}\right\}.
  \nonumber
  \end{eqnarray}
Let $M:=\left[2+\alpha+(\log T)/(\log 2)\right]$. In the following, we
suppose that $H\ge 1/2$ and $0<Z\le 2^{M}/T\le 2^{2+\alpha}$. By
$G(H,Z)$ we denote
the number of solutions to 

  \begin{displaymath}
  \left\vert \frac{k_1}{k_2}-
  \left(\frac{l_1+\theta}{l_2+\theta}\right)^{\alpha} \right\vert \le Z
  \end{displaymath}
with $k_1,k_2 \sim H$, $(k_1,k_2)=1$ and $l_1,l_2\sim L$. 
We then have
  \begin{eqnarray}
  & & \hspace{-0.4cm} 
  \sum\limits_{\scriptsize \begin{array}{cccc} &k_1,k_2\sim H,&\\ 
  &(k_1,k_2)=1& \end{array}} \hspace{-0.4cm} \sum\limits_{l_1,l_2\sim L}
  \min\left\{T,\ \left\vert \frac{k_1}{k_2}-
  \left(\frac{l_1+\theta}{l_2+\theta}\right)^{\alpha} 
  \right\vert^{-1}\right\}\label{61}\\ &\ll&
  T \sum\limits_{m=0}^{M} 
  G(H,2^m/T)2^{-m}.\nonumber
  \end{eqnarray}

Let $S(H)$ be the set of all fractions
$k_1/k_2$ with $k_1,k_2\sim H$, $(k_1,k_2)=1$. This set is well-spaced
with spacing $1/(4H^2)$. Hence,   
  \begin{eqnarray}
  G(H,Z)
  &= & \sum\limits_{l_1,l_2\sim L} 
  \left\vert \left\{u\in S(H) \ : \ 
  \left\vert
  \left(\frac{l_1+\theta}{l_2+\theta}\right)^{\alpha} - 
  u \right\vert \le Z\right\}
  \right\vert \label{62} \\ \nonumber\\ &\ll&
  L^2(ZH^2+1).\nonumber
  \end{eqnarray}
By a short calculation, from (\ref{62}), we derive
 
 \begin{equation}
  T \sum\limits_{m=0}^{M} 
  G(H,2^m/T)2^{-m} \ll (H^2L^2+TL^2)\log(2T).\label{63}
  \end{equation} 

We now estimate the left side of (\ref{63}) in an alternative way. 
Using the Cauchy-Schwarz inequality and taking the
above-mentioned spacing properties of the set $S(H)$ into account, we obtain
  \begin{eqnarray}
  & & G(H,Z) \label{64}\\ 
  &=&\sum\limits_{u\in S(H)} 
  \left\vert \left\{l_1,l_2\sim L\ : \ 
  \left\vert u-
  \left(\frac{l_1+\theta}{l_2+\theta}\right)^{\alpha} \right\vert \le Z\right\}
  \right\vert \nonumber\\ \nonumber\\ &\ll&
  H \left(\sum\limits_{u\in S(H)} 
  \left\vert \left\{l_1,l_2\sim L\ : \ 
  \left\vert u-
  \left(\frac{l_1+\theta}{l_2+\theta}\right)^{\alpha} \right\vert \le Z\right\}
  \right\vert^2\right)^{1/2} \nonumber\\ \nonumber\\
  &\ll& H\left((ZH^2+1)
  \left\vert \left\{l_1,l_2,l_1^{\prime},l_2^{\prime} \sim L\ : \ 
  \left\vert \left(\frac{l_1+\theta}{l_2+\theta}\right)^{\alpha}-\right.
  \right.\right.\right.\nonumber\\
  & & \left.\left.\left.\left. 
  \left(\frac{l_1^{\prime}+\theta}{l_2^{\prime}+\theta}\right)^{\alpha} 
  \right\vert
  \le 2Z\right\}
  \right\vert\right)^{1/2}. \nonumber
  \end{eqnarray}
Using Taylor$^{\prime}$s 
formula and $1/2\le (l_1+\theta)/(l_2+\theta)\le 2$, $1/2\le 
(l_1^{\prime}+\theta)/(l_2^{\prime}+\theta)\le 2$, we deduce that 
there is a positive 
constant $c$
depending only on $\alpha$ such
that
  \begin{eqnarray}
  & & \left\vert \left\{l_1,l_2,l_1^{\prime},l_2^{\prime} \sim L\ : \ 
  \left\vert \left(\frac{l_1+\theta}{l_2+\theta}\right)^{\alpha}-
  \left(\frac{l_1^{\prime}+\theta}{l_2^{\prime}+\theta}\right)^{\alpha} 
  \right\vert
  \le 2Z\right\}
  \right\vert \label{65} \\ \nonumber\\ 
  &\le& \left\vert \left\{l_1,l_2,l_1^{\prime},l_2^{\prime} \sim L\ : \ 
  \left\vert \frac{l_1+\theta}{l_2+\theta}-
  \frac{l_1^{\prime}+\theta}{l_2^{\prime}+\theta}\right\vert
  \le cZ\right\}
  \right\vert. \nonumber
  \end{eqnarray}
Taking the first inequality on page 145 of [Ha2] into account, we get
  \begin{eqnarray}
  & &\left\vert \left\{l_1,l_2,l_1^{\prime},l_2^{\prime} \sim L\ : \ 
  \left\vert \frac{l_1+\theta}{l_2+\theta}-
  \frac{l_1^{\prime}+\theta}{l_2^{\prime}+\theta}\right\vert
  \le cZ\right\} \right\vert \label{66}\\ \nonumber\\
  &\le & \vert \{l_1,l_2,l_1^{\prime},l_2^{\prime} \sim L\ : \ 
  \vert (l_1+\theta)(l_2^{\prime}+\theta)-
  (l_1^{\prime}+\theta)(l_2+\theta)\vert\nonumber\\
  & & \le cZ(2L+\theta)^2\} \vert 
  \nonumber\\ 
  &\ll& (ZL^2+1)L^2\log^2(2L). \nonumber
  \end{eqnarray}
The implied $\ll$-constant does not depend on $\theta$.
Combining (\ref{64}), (\ref{65}) and (\ref{66}), we get
 
 \begin{equation}
  G(H,Z)\ll (ZH^2L^2+Z^{1/2}(H+L)HL+HL)\log(2L). \label{67}
  \end{equation}
By a short calculation, from (\ref{67}), we derive
  \begin{eqnarray}
  & & T \sum\limits_{m=0}^{M} 
  G(H,2^m/T)2^{-m} \label{68} \\ &\ll& 
  (H^2L^2+T^{1/2}(H+L)HL+THL)\log(2L)\log(2T).\nonumber
  \end{eqnarray} 

Combining (\ref{63}) and (\ref{68}), and taking the condition $L\le T^{1/2}$ in
{\sc Theorem 4} into
account, we get
  \begin{eqnarray}
  & & T \sum\limits_{m=0}^{M} 
  G(H,2^m/T)2^{-m} \label{69}\\ &\ll&
  \left(H^2L^2+THL+T\min\left\{H^2,L^2\right\}\right)
  \log(2L)\log(2T).\nonumber
  \end{eqnarray} 
From (\ref{60}), (\ref{61}) and (\ref{69}), we obtain
 \begin{eqnarray}
  & & \label{70}\\
  & & I\ll \left(K^2L^2+TKL\log(2K)+T\sum\limits_{d\le 2K}
  \min\left\{\frac{K^2}{d^2},L^2\right\}\right)
  \log(2L)\log(2T). \nonumber 
  \end{eqnarray}
If $K\le L$, then we have

  \begin{displaymath}
  \sum\limits_{d\le 2K}
  \min\left\{\frac{K^2}{d^2},L^2\right\} \le \sum\limits_{d\le 2K}
  \frac{K^2}{d^2} \ll K^2\le KL.
  \end{displaymath}
Otherwise, we have

  \begin{displaymath}
  \sum\limits_{d\le 2K}
  \min\left\{\frac{K^2}{d^2},L^2\right\} \le \sum\limits_{d\le K/L} L^2 +
  \sum\limits_{K/L<d\le 2K} \frac{K^2}{d^2} \ll KL.
  \end{displaymath}
Therefore, from (\ref{70}) follows 

   \begin{displaymath}
  I\ll \left(K^2L^2+TKL\right)\log(2K)\log(2L)\log(2T).  
  \end{displaymath}
This implies the result of {\sc Theorem 4}. 
$\Box$\\

We note that if the condition $L\le T^{1/2}$ in
{\sc Theorem 4} could be removed 
for all $\theta\not=0$, then the condition
$\lambda\le 5/19$ in {\sc Theorem 3} could be removed 
for all irrational $\theta$.

\end{document}